\documentclass[11pt]{article}

\usepackage[T1]{fontenc}
\usepackage{lmodern}
\usepackage{amsmath,amssymb,amsthm,mathtools}
\usepackage[margin=1.08in]{geometry}
\usepackage{microtype}
\usepackage{enumitem}
\usepackage{booktabs}
\usepackage{xcolor}
\usepackage{hyperref}
\usepackage{fancyhdr}

\hypersetup{
  colorlinks=true,
  linkcolor=blue!55!black,
  citecolor=blue!55!black,
  urlcolor=blue!55!black,
  pdftitle={Mesh-Degree Rigidity for Positive Chebyshev-Fourier Approximants},
  pdfauthor={Vasily Stodolsky; Independent Researcher; ORCID 0009-0003-5329-5308},
  pdfsubject={Positive Chebyshev-Fourier approximation; AI-assisted research disclosed in full text},
  pdfkeywords={entire functions, trigonometric polynomials, Chebyshev polynomials, Laguerre-Polya class, weak convergence, exponential type}
}

\newtheorem{theorem}{Theorem}[section]
\newtheorem{proposition}[theorem]{Proposition}
\newtheorem{corollary}[theorem]{Corollary}
\newtheorem{remark}[theorem]{Remark}
\newtheorem{example}[theorem]{Example}

\newcommand{\R}{\mathbb{R}}
\newcommand{\C}{\mathbb{C}}
\newcommand{\dd}{\,\mathrm{d}}
\newcommand{\supp}{\operatorname{supp}}
\newcommand{\ord}{\operatorname{ord}}
\newcommand{\e}{\mathrm{e}}
\newcommand{\weakto}{\Rightarrow}
\newcommand{\T}{\mathsf{T}}

\title{\textbf{Mesh-Degree Rigidity for Positive\\
  Chebyshev-Fourier Approximants}}
\author{Vasily Stodolsky\\
\small Independent Researcher\\
\small \href{https://orcid.org/0009-0003-5329-5308}
{ORCID 0009-0003-5329-5308}}
\date{Version 0.1.4, August 7, 2026\\arXiv-preparation and provenance update}

\begin{document}

\maketitle

\begin{center}
\begin{minipage}{0.88\textwidth}
\small
\textbf{Authorship and verification status.}
Vasily Stodolsky is the sole author of this manuscript. AI systems were used as
research, writing, and verification tools and are not authors. Vasily Stodolsky
directed the workflow, set its scope and release criteria, curated the
mathematical and verification record, made the final decisions on claims and
revisions, approved this version, and accepts responsibility for its contents.
The manuscript contains complete theorem statements, candidate proofs, and
boundary examples developed in a materially AI-assisted workflow. Detailed
model roles and verification limitations are reported in the end-matter
disclosure. No independent human peer review is claimed.

\smallskip
\textbf{Version note.}
Version 0.1.4 adds a nontechnical proof guide and an independently
model-rechecked non-vacuity construction, clarifies human responsibility, and
records a background attribution. The theorem statements and formal proofs are
unchanged from version 0.1.3.
\end{minipage}
\end{center}
\medskip

\begin{abstract}
Let $d_k>0$ and let $Q_k$ be a polynomial of degree $r_k$ with
nonnegative Chebyshev coefficients.  This artifact studies locally uniform limits of
\[
  f_k(z)=C_k Q_k(\cos(d_k z)), \qquad C_k>0.
\]
For a locally uniform limit $F$ with $F(0)>0$, the first result records the
measure-theoretic rigidity supplied by complex compact-open convergence: the positive lattice measures
represented by the $f_k$ converge weakly, their second moments converge, and
their quadratic tails are uniformly integrable.  If $f_k\to F$ locally
uniformly with $F(0)>0$, all zeros of $Q_k$ in $[-1,1)$, $d_k\to0$, and
$\ord(F)<2$, then the main estimate couples the frequency mesh, the
polynomial degree, and the nonzero real zeros $\{\pm\gamma_n\}$ of $F$,
counted with multiplicity:
\[
  \limsup_{k\to\infty} r_kd_k^2
  \leq 4\sum_{\gamma_n>A}\frac{1}{\gamma_n^2}
\]
for every $A>0$ with $F(A)F(-A)\neq0$.  Combining this estimate with the
period and positivity arguments, if $f_k\to F$ locally uniformly with
$F(0)>0$, all zeros of $Q_k$ lie in $[-1,1)$, and $F$ has no nonzero real
period, has order below two, and is not of finite exponential type, then
\[
  d_k\longrightarrow0,\qquad
  r_kd_k\longrightarrow\infty,\qquad
  r_kd_k^2\longrightarrow0.
\]
Equivalently, the degree is forced into the window
$d_k^{-1}=o(r_k)$ and $r_k=o(d_k^{-2})$.  Elementary examples exhibit the
role of each hypothesis and the Gaussian boundary at order two.
\end{abstract}



\medskip
\noindent\textbf{Verification materials and contact.}
The audit summary is recorded in the release metadata and workflow
disclosure.  The release package also contains the executable checks, source,
and SHA-256 manifest.  The persistent version history is maintained at
\href{https://doi.org/10.5281/zenodo.21672811}
{10.5281/\allowbreak zenodo.21672811}.
Corrections may be reported to Vasily Stodolsky via
\href{https://x.com/Vasya947}{@Vasya947}.

\medskip
\noindent\textbf{License.}
The manuscript and documentation are licensed under
\href{https://creativecommons.org/licenses/by/4.0/}{CC BY 4.0}; executable
code in the release is licensed under the MIT License.

\section{Introduction}

Compact-open limits of polynomials with restricted zeros are governed by
classical canonical-product theory; see, for example, Levin
\cite[Chapter VIII, Section 1, Theorem 3]{Levin1980}.  The approximation
results of P\'olya \cite{Polya1913},
Lindwart-P\'olya \cite{LindwartPolya1914}, and the work of Korevaar
\cite{Korevaar1951} provide the basic framework.  In modern language, a
locally uniform limit of real-rooted polynomials belongs to the
Laguerre-P\'olya class and may acquire a zero-free Gaussian factor from roots
that escape to infinity; see also \cite{Karlin1968}.
Recent formulations make the roles of the limiting reciprocal roots and the
first and second power sums of reciprocal roots especially explicit
\cite{Assiotis2022}.

There is a related fixed-period theory for trigonometric polynomials.
Prather \cite[Theorem 2.1]{Prather1979} proved that a finite-order compact
limit of balanced $2\pi$-periodic trigonometric polynomials with only real
zeros is itself a trigonometric polynomial.  The regime considered here is
different: the period $2\pi/d_k$ varies with $k$ and, for a nonperiodic
limit, must tend to infinity.

The purpose of this note is not to introduce a new canonical
factorization.  It is to extract a quantitative consequence of the special
composition
\[
  Q_k(\cos(d_kz)).
\]
Each zero of $Q_k$ in $[-1,1)$ contributes at least $d_k^2/2$ to a finite
logarithmic-derivative identity.  Zeros that remain in a fixed $z$-window
converge, with multiplicity, to zeros of the limiting entire function.
The remaining contribution is therefore bounded by the reciprocal-square
tail of the limiting zero divisor.  If all zeros of $Q_k$ lie in $[-1,1)$,
$f_k\to F$ locally uniformly with $F(0)>0$, $d_k\to0$, and $\ord(F)<2$,
this gives
\[
  \limsup_{k\to\infty}r_kd_k^2
  \leq 4\sum_{\gamma_n>A}\gamma_n^{-2}
\]
for every $A>0$ with $F(A)F(-A)\neq0$.  The right-hand side tends to zero
as $A\to\infty$ through such values.

Nonnegative Chebyshev coefficients supply a second, independent mechanism.
They identify $f_k$ as the Fourier transform of a positive lattice measure
supported in $[-r_kd_k,r_kd_k]$.  If $r_kd_k$ stayed bounded, the limit
would have finite exponential type.  Combining the two mechanisms produces
the forced degree window
\[
  d_k^{-1}=o(r_k),\qquad r_k=o(d_k^{-2}).
\]

The measure argument is useful independently of the zero geometry.  We
therefore first state it for arbitrary positive even compactly supported
measures.  It is a direct consequence of the continuity theorem, complex
compact-open convergence, and convergence of derivatives; it is included
to make the precise moment and interval conclusions transparent.  Related
factorization for Lee-Yang variables appears already in
\cite[Proposition 2]{Newman1975}, and a Fourier-transform treatment is given
by Cardon \cite{Cardon2005}.  Weak-limit closure for Lee-Yang variables is
given in \cite[Theorem 7 in arXiv version 3]{NewmanWu2019}; a locally uniform Hurwitz step already
appears in \cite[proof of Theorem 10, p.~7]{Newman1975}.  The classical
continuity theorem for moment-generating functions goes back at least to
Curtiss \cite{Curtiss1942}.  The fixed-radius scaling $r_kd_k=O(1)$ also
belongs to the broader landscape of classical Mehler-Heine asymptotics; see
\cite[Theorem 8.1.1]{Szego1975}.  Example \ref{ex:dirichlet} below is obtained
directly from the Dirichlet-kernel identity.

For background on related Laguerre-P\'olya and probabilistic representations,
see also Konstantopoulos, Patie, and Sarkar
\cite[Theorem~3.14]{KonstantopoulosPatieSarkar2024}, \emph{Annales de l'Institut
Fourier} \textbf{74}(1) (2024), 377--421, DOI
\href{https://doi.org/10.5802/aif.3600}{10.5802/aif.3600}. This reference is
included as contextual attribution only; no result from it is used as a step
in the formal proofs below.

\section*{Proof Idea}

\textbf{Purpose of this guide.}
This section is an explanatory map to the already stated results. It does not
replace a formal proof or enlarge the scope of any hypothesis. In addition to
the map, it records a separately labelled non-vacuity construction outside the
frozen formal region; that construction was accepted by model recheck and is
released for independent verification. The formal definitions, statements,
and proofs begin in the next section. The guide is included so that a reader
can distinguish the two mechanisms used by the paper and can see where each
assumption enters.

\textbf{Starting point.}
The composition $Q_k(\cos(d_k z))$ carries two kinds of structure at the
same time. The condition on the Chebyshev coefficients says that, after the
normalization factor $C_k$, the function is the Fourier transform of a
positive even measure supported on the lattice
$\{md_k:-r_k\leq m\leq r_k\}$. The condition that the zeros of $Q_k$ lie
in $[-1,1)$ says that the zeros of the composed entire function occur in
controlled periodic families. Compact-open convergence transfers information
from both descriptions to the limit $F$.

\textbf{First mechanism: the positive measure.}
Nonnegative Chebyshev coefficients give a positive measure $\nu_k$ whose
Fourier transform is $f_k$. Compact-open convergence is stronger than merely
pointwise convergence on the real line: it also controls derivatives at the
origin. The measure part of the proof uses this to obtain weak convergence of
the measures, convergence of the relevant second moments, and control of
their quadratic tails. Positivity is essential here. Without it, cancellation
can make a Fourier transform look controlled while the underlying signed
mass behaves very differently.

The support of $\nu_k$ is contained in
$[-r_kd_k,r_kd_k]$. If those support radii remained bounded, the limiting
Fourier transform would have finite exponential type. Thus the hypothesis
that $F$ is not of finite exponential type forces the support radius, and
hence $r_kd_k$, to escape to infinity. This is the lower side of the degree
window. It is a support-versus-type argument, not a count of zeros.

\textbf{Second mechanism: the zero budget.}
The zero restriction on $Q_k$ gives a logarithmic-derivative identity for
the composition. Each root of $Q_k$ contributes a nonnegative amount, with a
uniform contribution of order $d_k^2$ after the relevant normalization. The
total is therefore a budget that can be compared with the zero divisor of
the compact-open limit. Zeros that stay in a fixed bounded window are matched
with zeros of $F$, including multiplicity. What is left after the fixed
window is bounded by the reciprocal-square tail of the zeros of $F$.

This is why the estimate is written with a cutoff $A$. For a fixed cutoff,
one first passes to the limit in $k$ and obtains a bound for
$r_kd_k^2$ by a tail of the limit divisor. Only after that may the cutoff be
sent outward. The assumption that $F$ has order below two makes that
reciprocal-square tail tend to zero along admissible cutoffs. This produces
the upper side of the degree window, namely $r_kd_k^2\to0$.

\textbf{Logical route and hypotheses.}
The no-nonzero-period assumption rules out a fixed positive mesh scale: a
subsequential positive limit for $d_k$ would retain a period in the limit.
Hence $d_k\to0$. The non-finite-type assumption then turns the positive
measure representation into $r_kd_k\to\infty$. Finally, the order-below-two
assumption turns the root-budget tail into $r_kd_k^2\to0$. Together these are
exactly the two asymptotic inequalities expressed by
$d_k^{-1}=o(r_k)$ and $r_k=o(d_k^{-2})$.

\textbf{What this guide does not claim.}
It does not construct an approximating sequence for a chosen $F$, establish
a rate of convergence, identify individual roots outside fixed windows, or
permit a diagonal cutoff $A=A(k)$. The following construction addresses the
non-vacuity of the joint hypotheses. It was accepted by an independent focused
model recheck, but it is not a general existence theorem for prescribed targets.

\textbf{Non-vacuity construction.}
For $k\geq1$, set $d_k=1/k$,
\[
  m_{j,k}=\lfloor k j^{-2/3}\rfloor,\qquad
  Q_k(x)=\prod_{j=1}^k T_{m_{j,k}}(x),\qquad C_k=1,
\]
and write $r_k=\sum_{j=1}^k m_{j,k}$. Then
\[
  Q_k(\cos(d_kz))=\prod_{j=1}^k
  \cos\!\left(\frac{m_{j,k}}{k}z\right),
  \qquad
  F(z)=\prod_{j=1}^{\infty}\cos(j^{-2/3}z).
\]
Every $m_{j,k}$ is positive. The roots of each Chebyshev factor lie in
$(-1,1)$, and the identity
$T_aT_b=(T_{a+b}+T_{|a-b|})/2$ shows inductively that $Q_k$ has nonnegative
Chebyshev coefficients. For each fixed factor index $j$,
$m_{j,k}/k\to j^{-2/3}$. On a fixed disk, first compare any finite initial
product factor by factor, and then choose that finite cutoff so that the
remaining square-frequency tail $\sum_{j>J}j^{-4/3}$ is small. The rounding
error in the initial product tends to zero as $k\to\infty$. This fixed-factor,
then-tail route gives locally uniform convergence to $F$ by compact-product
estimates; it does not use an unjustified diagonal cutoff.

The zeros of the factors give a zero-counting function of order
$N_F(R)\asymp R^{3/2}$. Splitting the product at
$j\asymp R^{3/2}$ also gives $\log M_F(R)=O(R^{3/2})$, so $F$ has order
exactly $3/2$. A nonzero real period would make the real zero divisor repeat
from one compact period interval and would force $N_F(R)=O(R)$, a
contradiction. Since the order exceeds one, $F$ is not of finite exponential
type. Finally, integral comparison gives $r_k\asymp k^{4/3}$, hence
$r_kd_k\to\infty$ and $r_kd_k^2\to0$. An independent focused GPT-5.5 xhigh
mathematical recheck accepted each hypothesis of Theorem~\ref{thm:window} for
this construction.

\section{Setup and main results}

\subsection{Normalization}

For each $k$, let $d_k>0$, $C_k>0$, and let $Q_k$ be a real polynomial of
exact degree $r_k\geq1$.  We use the normalization
\begin{equation}\label{eq:cheb-expansion}
  Q_k(x)
  =a_{0,k}+2\sum_{m=1}^{r_k}a_{m,k}\T_m(x),
  \qquad
  a_{m,k}\geq0,\quad a_{r_k,k}>0,
\end{equation}
where $\T_m(\cos\theta)=\cos(m\theta)$.  Set
\begin{equation}\label{eq:fk}
  f_k(z)=C_kQ_k(\cos(d_kz))
\end{equation}
and define the finite positive even measure
\begin{equation}\label{eq:nuk}
  \nu_k=C_k\left[
    a_{0,k}\delta_0+
    \sum_{m=1}^{r_k}a_{m,k}
    \bigl(\delta_{md_k}+\delta_{-md_k}\bigr)
  \right].
\end{equation}
With the Fourier convention
\[
  \widehat{\nu}(z)=\int_{\R}\e^{izu}\dd\nu(u),
\]
we have
\begin{equation}\label{eq:dictionary}
  f_k=\widehat{\nu_k},\qquad
  \nu_k(\R)=f_k(0),\qquad
  \supp\nu_k\subset[-r_kd_k,r_kd_k],\qquad
  \int_{\R}u^2\dd\nu_k(u)=-f_k''(0).
\end{equation}
Since $a_{r_k,k}>0$, the extreme atoms at $\pm r_kd_k$ in \eqref{eq:nuk}
show that $r_kd_k$ is the actual support radius.  We also call it the
frequency radius, while $d_k$ is the frequency mesh.

Convergence of entire functions will always mean locally uniform convergence
on $\C$.  For an entire function $F$, write
\[
  M_F(R)=\max_{|z|=R}|F(z)|,\qquad
  \ord(F)=\limsup_{R\to\infty}
  \frac{\log\log M_F(R)}{\log R}.
\]
We say that $F$ has finite exponential type if there are constants
$A,B<\infty$ such that
\[
  |F(z)|\leq A\e^{B|z|},\qquad z\in\C.
\]

\subsection{Positive-measure rigidity}

\begin{proposition}[Measure and moment rigidity]\label{prop:measure}
Let $\mu_k$ be finite positive even compactly supported measures on $\R$ and
let
\[
  g_k(z)=\int_{\R}\e^{izu}\dd\mu_k(u).
\]
Suppose that $g_k\to G$ locally uniformly on $\C$ and $G(0)>0$.  Then there
is a unique finite positive even measure $\mu$ such that
\[
  G(z)=\int_{\R}\e^{izu}\dd\mu(u),\qquad z\in\C.
\]
Moreover:
\begin{enumerate}[label=\textup{(\roman*)}]
  \item $\mu_k\weakto\mu$ weakly as finite measures;
  \item $\mu_k(B)\to\mu(B)$ for every Borel set $B$ satisfying
        $\mu(\partial B)=0$;
  \item
  \[
    \int_{\R}u^2\dd\mu_k(u)
    \longrightarrow
    \int_{\R}u^2\dd\mu(u)=-G''(0);
  \]
  \item the quadratic tails are uniformly integrable:
  \[
    \lim_{A\to\infty}\sup_k
    \int_{|u|>A}u^2\dd\mu_k(u)=0;
  \]
  \item for every continuous $h$ with
        $|h(u)|\leq C_h(1+u^2)$,
  \[
    \int_{\R}h\dd\mu_k\longrightarrow\int_{\R}h\dd\mu;
  \]
  \item
  \[
    \lim_{A\to\infty}\sup_k\mu_k(\{|u|>A\})=0.
  \]
\end{enumerate}
The limit $\mu$ has all exponential moments:
\[
  \int_{\R}\e^{R|u|}\dd\mu(u)<\infty
  \qquad\text{for every }R>0.
\]
\end{proposition}

\begin{corollary}[Lattice masses]\label{cor:lattice}
Suppose that \eqref{eq:cheb-expansion} through \eqref{eq:nuk} hold and that
$f_k\to F$ locally uniformly with $F(0)>0$.  Let $\nu$ be the measure from
Proposition \ref{prop:measure}.  If $I\subset(0,\infty)$ is an interval and
$\nu(\partial I)=0$, then
\[
  C_k\sum_{md_k\in I}a_{m,k}\longrightarrow\nu(I).
\]
In addition,
\begin{align*}
  \lim_{A\to\infty}\sup_k
  C_k\sum_{md_k>A}a_{m,k}&=0,\\
  \lim_{A\to\infty}\sup_k
  C_k\sum_{md_k>A}(md_k)^2a_{m,k}&=0.
\end{align*}
If $\nu$ is atomless, the conclusion holds for every interval
$I\subset(0,\infty)$.
\end{corollary}

\subsection{The root budget and the forced window}

Assume from now on that every zero of $Q_k$, counted with multiplicity,
belongs to $[-1,1)$.  It follows that all zeros of $f_k$ are real.

\begin{theorem}[Quantitative root budget]\label{thm:budget}
Suppose that $f_k\to F$ locally uniformly, $F(0)>0$, $d_k\to0$, and
$\ord(F)<2$.  List the positive zeros of $F$, with multiplicity and one
representative from each pair, as
\[
  0<\gamma_1\leq\gamma_2\leq\cdots .
\]
The list is allowed to be finite or empty.  Then
\[
  \sum_n\frac{1}{\gamma_n^2}<\infty,
\]
and for every $A>0$ such that $F(A)F(-A)\neq0$,
\begin{equation}\label{eq:budget}
  \limsup_{k\to\infty}r_kd_k^2
  \leq4\sum_{\gamma_n>A}\frac{1}{\gamma_n^2}.
\end{equation}
In particular,
\[
  r_kd_k^2\longrightarrow0.
\]
\end{theorem}

\begin{remark}
The proof of Theorem \ref{thm:budget} uses the location of the zeros of
$Q_k$, but not the nonnegativity of the coefficients in
\eqref{eq:cheb-expansion}.  Positivity enters the lower scale law below.
\end{remark}

\begin{remark}[A classical subcase]
If the stronger root restriction $x_{j,k}\in[-1,0]$ holds, the qualitative
conclusion $r_kd_k^2\to0$ also follows from classical triangular-array
probability.  Indeed, every normalized factor
\[
  \frac{\cos(d_kz)-x_{j,k}}{1-x_{j,k}}
\]
is the characteristic function of a variable supported on
$\{0,\pm d_k\}$ with variance $d_k^2/(1-x_{j,k})\in[d_k^2/2,d_k^2]$.
The array is infinitesimal.  The classical convergence criterion
\cite[Chapter V, Section 26, Theorem 1, p.~126]{GnedenkoKolmogorov1954}
then makes any weak limit normal; the strict order bound excludes a
nondegenerate Gaussian limit.  Hence the total variance, and therefore
$r_kd_k^2$, tends to zero.  This argument does not cover roots in $(0,1)$,
where the individual normalized factors correspond to signed rather than
positive measures.
\end{remark}

\begin{theorem}[Forced mesh-degree window]\label{thm:window}
In addition to the standing assumptions of Sections 2.1 and 2.3, suppose that
$f_k$ has the form \eqref{eq:fk}, with $Q_k$ normalized as in
\eqref{eq:cheb-expansion}, with $\nu_k$ as in \eqref{eq:nuk}, and with all
zeros of $Q_k$ in $[-1,1)$.  Assume further that:
\begin{enumerate}[label=\textup{(\roman*)}]
  \item $f_k\to F$ locally uniformly and $F(0)>0$;
  \item $F$ has no nonzero real period;
  \item $\ord(F)<2$;
  \item $F$ is not of finite exponential type.
\end{enumerate}
Then
\begin{equation}
  d_k\longrightarrow0,\qquad
  r_kd_k\longrightarrow\infty,\qquad
  r_kd_k^2\longrightarrow0.
\end{equation}
Equivalently,
\[
  d_k^{-1}=o(r_k),\qquad r_k=o(d_k^{-2}).
\]
\end{theorem}

\section{Proof of the measure proposition}

Put $M_k=\mu_k(\R)=g_k(0)$ and $M=G(0)$.  Then $M_k\to M>0$.
For all sufficiently large $k$, let $p_k=\mu_k/M_k$.  On the real axis,
their characteristic functions converge pointwise to
\[
  \phi(t)=\frac{G(t)}{M},
  \qquad \phi(0)=1.
\]
Each $g_k$ is even, so local uniform convergence makes $G$, and hence
$\phi$, even.  The function $\phi$ is continuous at zero.  The continuity theorem for
characteristic functions therefore gives a probability measure $p$ such
that $p_k\weakto p$ and $\widehat p=\phi$ on $\R$; see, for example,
\cite[Theorem 26.3 and the corollary following it]{Billingsley1995}.  Set
$\mu=Mp$.  This proves weak
convergence of the finite measures.  Since $\phi$ is even, the reflection
of $p$ has the same characteristic function as $p$; uniqueness therefore
shows that $p$, and hence $\mu$, is even.

Fix $y\in\R$.  Evenness of the measures identifies the Fourier transform on
the imaginary axis with the integral of $\cosh(yu)$.  Since
$u\mapsto\cosh(yu)$ is nonnegative and lower semicontinuous, the Portmanteau
theorem gives
\begin{align*}
  \int_{\R}\cosh(yu)\dd\mu(u)
  &\leq\liminf_{k\to\infty}
      \int_{\R}\cosh(yu)\dd\mu_k(u)\\
  &=\liminf_{k\to\infty}g_k(iy)
   =G(iy)<\infty.
\end{align*}
Alternatively, one may first apply weak convergence to the bounded continuous
truncations $\min\{\cosh(yu),N\}$ and then use monotone convergence.  Since
$\e^{R|u|}\leq2\cosh(Ru)$, the measure $\mu$ has every exponential moment.
Consequently
\[
  H(z)=\int_{\R}\e^{izu}\dd\mu(u)
\]
is entire.  Indeed, on $|z|\leq R$, every differentiated integrand is
dominated by an integrable function because
\[
  |u|^n\e^{R|u|}
  \leq c_n\e^{(R+1)|u|}.
\]
On the real axis $H=G$, and hence the identity theorem yields $H=G$ on
$\C$.  Uniqueness follows from uniqueness of a finite measure determined
by its characteristic function.

Since locally uniform convergence of holomorphic functions implies
convergence of every derivative,
\[
  \int_{\R}u^2\dd\mu_k(u)
  =-g_k''(0)
  \longrightarrow-G''(0)
  =\int_{\R}u^2\dd\mu(u).
\]
It remains to record why this entails uniform integrability.  For $R>0$,
let $h_R(u)=\min(u^2,R^2)$.  Since
$(u^2-R^2)_+=u^2-h_R(u)$, weak convergence and convergence of the full
second moments give
\[
  \int_{\R}(u^2-R^2)_+\dd\mu_k(u)
  \longrightarrow
  \int_{\R}(u^2-R^2)_+\dd\mu(u).
\]
For $|u|>2R$,
\[
  u^2\leq\frac43(u^2-R^2).
\]
It follows that
\[
  \limsup_{k\to\infty}
  \int_{|u|>2R}u^2\dd\mu_k(u)
  \leq\frac43
  \int_{\R}(u^2-R^2)_+\dd\mu(u),
\]
and the right-hand side tends to zero as $R\to\infty$.  The finitely many
initial indices are handled separately.  This proves the uniform
integrability assertion.

For a continuous $h$ of quadratic growth, multiply by a continuous cutoff
that is one on $[-R,R]$ and zero outside $[-2R,2R]$.  Weak convergence
handles the cutoff part, while
\[
  \mu_k(\{|u|>R\})
  \leq R^{-2}\int_{|u|>R}u^2\dd\mu_k(u)
\]
and uniform integrability handle the tails.  The same estimate gives
uniform tightness.  Finally, convergence on continuity sets follows from
Portmanteau.  Proposition \ref{prop:measure} is proved.  Corollary
\ref{cor:lattice} follows by applying Proposition \ref{prop:measure}(ii),
(iv), and (vi) to the measures \eqref{eq:nuk}.

\section{Proof of the scale theorems}

Sections 4.2 and 4.3 establish two conclusions needed for Theorem
\ref{thm:window}.  For Theorem \ref{thm:budget}, $d_k\to0$ is already a
hypothesis; Section 4.5 uses that hypothesis directly.  Sections 4.4 and 4.5
do not use the nonperiodicity or infinite-type assumptions.

\subsection{A finite logarithmic-derivative identity}

Suppress the index $k$.  Write the zeros of $Q$, with multiplicity, as
\[
  x_j=\cos\alpha_j,\qquad
  \alpha_j\in(0,\pi],\qquad 1\leq j\leq r,
\]
and put $t_j=\alpha_j/d$.  Since $Q(1)\neq0$,
\[
  \frac{Q'(1)}{Q(1)}
  =\sum_{j=1}^r\frac{1}{1-x_j}.
\]
Differentiating $f(z)=C Q(\cos(dz))$ twice at zero gives
$f''(0)=-Cd^2Q'(1)$.  Therefore
\begin{equation}\label{eq:finite-identity}
  -\frac{f''(0)}{f(0)}
  =\sum_{j=1}^r
  \frac{d^2}{2\sin^2(dt_j/2)}.
\end{equation}
Every summand is at least $d^2/2$.  If $x_j=-1$, then $t_j=\pi/d$ and
the corresponding summand is exactly $d^2/2$.  Although the resulting
zero of the composition is double, \eqref{eq:finite-identity} counts the
root of $Q$ once per algebraic multiplicity.

\subsection{The mesh tends to zero}

Every $f_k$ has the real period $L_k=2\pi/d_k$.  Suppose that $d_k$ does
not tend to zero.  Along a subsequence, $d_k\geq\varepsilon>0$, and after
passing to a further subsequence, $L_k\to L\in[0,\infty)$.

If $L>0$, local uniform convergence in
\[
  f_k(z+L_k)=f_k(z)
\]
gives $F(z+L)=F(z)$, a forbidden nonzero real period.  If $L=0$, then
\[
  0=\frac{f_k(z+L_k)-f_k(z)}{L_k}
   =\frac1{L_k}\int_0^{L_k}f_k'(z+s)\dd s.
\]
The derivatives converge locally uniformly, so $F'(z)=0$.  A constant
function has every real period, again a contradiction.  Thus
\begin{equation}\label{eq:d-zero}
  d_k\longrightarrow0.
\end{equation}

\subsection{The frequency radius diverges}

Let $U_k=r_kd_k$.  Positivity and \eqref{eq:dictionary} imply
\begin{equation}\label{eq:type-bound}
  |f_k(z)|
  \leq f_k(0)\e^{U_k|\operatorname{Im}z|}.
\end{equation}
If $U_k$ failed to tend to infinity, there would be a subsequence on which
$U_k\leq U$.  Passing to the limit in \eqref{eq:type-bound} would give
\[
  |F(z)|\leq F(0)\e^{U|\operatorname{Im}z|}
  \leq F(0)\e^{U|z|},
\]
so $F$ would have finite exponential type.  Hence
\begin{equation}\label{eq:radius-infty}
  r_kd_k\longrightarrow\infty.
\end{equation}

\subsection{The canonical product of the limit}

All zeros of the $f_k$ are real.  Hurwitz's theorem shows that all zeros of
the nonzero limit $F$ are real as well.  The limit is even, real on $\R$,
and nonzero at the origin.  Since $\ord(F)<2$, Hadamard factorization,
paired over the zeros $\pm\gamma_n$, yields
\begin{equation}\label{eq:canonical}
  F(z)=F(0)\prod_n\left(1-\frac{z^2}{\gamma_n^2}\right),
  \qquad
  \sum_n\frac{1}{\gamma_n^2}<\infty.
\end{equation}
The possible exponential remainder has degree at most one; evenness removes
its linear part.  A quadratic Gaussian factor is excluded by the strict
order bound.  Differentiating at zero gives
\begin{equation}\label{eq:limit-second}
  -\frac{F''(0)}{F(0)}
  =2\sum_n\frac{1}{\gamma_n^2}.
\end{equation}
The factorization with a possible nonnegative quadratic exponential term
appears in the Lee-Yang setting in \cite[Proposition 2]{Newman1975}.  The
strict-order implication that this Gaussian term vanishes is stated
explicitly in \cite[Lemma 3.4, p.~166]{LiebSokal1981}.  Compare also the
variance decomposition in
\cite[Proposition 13 in arXiv version 3]{NewmanWu2019}.

\subsection{A fixed window of zeros}

Fix $A>0$ such that $F(A)F(-A)\neq0$.  All zeros of $F$ are real, so the
two endpoint conditions imply that $F$ is zero-free on $|z|=A$.  Local
uniform convergence makes $f_k$ zero-free in a neighbourhood
of that circle for all large $k$; in particular, no $t_{j,k}$ equals $A$.
By the assumed or previously established convergence $d_k\to0$,
$A<\pi/d_k$ for all large $k$.  Inside the disk $|z|<A$, the zeros of
$f_k$ are precisely the pairs
\[
  \pm t_{j,k},\qquad
  t_{j,k}=\frac{\alpha_{j,k}}{d_k}<A.
\]
Neither a periodic copy nor the endpoint zero $\pi/d_k$ lies in the disk.

The argument principle on the circle $|z|=A$ shows that these zeros converge,
with multiplicity, to the zeros $\pm\gamma_n$ of $F$ in the disk.  Therefore
\begin{equation}\label{eq:near-window}
  \sum_{t_{j,k}<A}
  \frac{d_k^2}{2\sin^2(d_kt_{j,k}/2)}
  \longrightarrow
  2\sum_{\gamma_n<A}\frac{1}{\gamma_n^2}.
\end{equation}
On the other hand, convergence of derivatives and values at zero, together
with \eqref{eq:finite-identity} and \eqref{eq:limit-second}, gives convergence
of the full sums.  Subtracting \eqref{eq:near-window} yields
\begin{equation}\label{eq:far-window}
  \sum_{t_{j,k}>A}
  \frac{d_k^2}{2\sin^2(d_kt_{j,k}/2)}
  \longrightarrow
  2\sum_{\gamma_n>A}\frac{1}{\gamma_n^2}.
\end{equation}

Let $n_k(A)$ denote the number of roots with $t_{j,k}<A$.  This number is
eventually constant and in particular bounded.  Since $d_k\to0$, it follows
that $n_k(A)d_k^2\to0$.  Since every remaining
summand in \eqref{eq:far-window} is at least $d_k^2/2$,
\[
  \frac{(r_k-n_k(A))d_k^2}{2}
  \leq
  \sum_{t_{j,k}>A}
  \frac{d_k^2}{2\sin^2(d_kt_{j,k}/2)}.
\]
Taking the upper limit proves \eqref{eq:budget}.  Finally let
$A\to\infty$ through values that are not zeros of $F$.  The tail in
\eqref{eq:budget} tends to zero by \eqref{eq:canonical}, proving
Theorem \ref{thm:budget}.  Equations \eqref{eq:d-zero},
\eqref{eq:radius-infty}, and Theorem \ref{thm:budget} prove
Theorem \ref{thm:window}.

\section{Boundary examples}

The examples below also serve to distinguish the three mechanisms in
Theorem \ref{thm:window}.

\begin{example}[A fixed period]
Let $d_k=1$, $Q_k(x)=1+x$, and $C_k=1$.  Then
\[
  f_k(z)=1+\cos z
\]
is a nonconstant periodic limit, while $d_k$ does not tend to zero.  The
Chebyshev coefficients are nonnegative and the only root of $Q_k$ is $-1$.
\end{example}

\begin{example}[A bounded frequency radius]\label{ex:dirichlet}
Let
\[
  D_r(x)=1+2\sum_{m=1}^r\T_m(x),\qquad
  d_r=\frac1r,\qquad C_r=\frac1r.
\]
The roots of $D_r$ are
\[
  \cos\left(\frac{2\pi\ell}{2r+1}\right),
  \qquad 1\leq\ell\leq r,
\]
and belong to $(-1,1)$.  The Dirichlet-kernel identity gives
\[
  f_r(z)=\frac1r
  \frac{\sin((r+\tfrac12)z/r)}{\sin(z/(2r))}
  \longrightarrow\frac{2\sin z}{z}
\]
locally uniformly, with the value at zero understood by continuity.
Here $r d_r=1$.  The limit is nonperiodic and has order one, but it has
finite exponential type.  Thus the infinite-type assumption is essential
for the lower scale law.
\end{example}

\begin{example}[The Gaussian boundary]
Let
\[
  d_N=N^{-1/2},\qquad r_N=N,\qquad
  Q_N(x)=(2x)^N,\qquad C_N=2^{-N}.
\]
The Chebyshev coefficients of $(2x)^N$ are nonnegative, and all of its
zeros are at the origin.  Moreover,
\[
  f_N(z)=\cos(z/\sqrt N)^N
  \longrightarrow\e^{-z^2/2}
\]
locally uniformly.  Here
\[
  r_Nd_N\longrightarrow\infty,\qquad
  r_Nd_N^2=1.
\]
The limiting Gaussian factor has order two.  This is the sharp boundary for
the strict order assumption in the upper scale law.  The convergence is the
classical de Moivre-Laplace limit; the same approximation is recorded in the
ferromagnetic limit-distribution setting in
\cite[p.~146]{SimonGriffiths1973}.
\end{example}

\begin{example}[Invisible degree outside the root interval]
Let $D_N$ be as in Example \ref{ex:dirichlet}, and define
\[
  H_N(x)=
  \left(
    \frac{1+N^{-1}(1+x)}{1+2/N}
  \right)^{N^2}.
\]
Every zero of $H_N$ equals $-N-1$, outside $[-1,1)$.  The factor
$1+N^{-1}(1+x)$ has nonnegative coefficients in the basis
$\T_0,\T_1$.  Powers and products preserve nonnegative Chebyshev
coefficients because
\[
  \T_m(x)\T_n(x)
  =\frac12\bigl(\T_{m+n}(x)+\T_{|m-n|}(x)\bigr).
\]
For $d_N=1/N$,
\[
  H_N(\cos(z/N))\longrightarrow1
\]
locally uniformly.  Consequently
\[
  \frac1N D_N(\cos(z/N))H_N(\cos(z/N))
  \longrightarrow\frac{2\sin z}{z},
\]
whereas the degree is $N+N^2$ and
\[
  (N+N^2)d_N^2\longrightarrow1.
\]
Thus the root restriction in Theorem \ref{thm:budget} cannot simply be
discarded, even if coefficient positivity is retained.
\end{example}

\begin{remark}[No microscopic root rigidity]
The conclusions do not locate individual roots.  Take the Dirichlet
approximant in Example \ref{ex:dirichlet} with $r=N$.  Multiplying the
$z$-domain expression by
\[
  \cos(z/N)^{s_N},\qquad s_N=\lfloor N^{3/2}\rfloor,
\]
corresponds to multiplying the underlying polynomial by $x^{s_N}$.  It adds
$s_N$ roots at $x=0$, preserves nonnegative Chebyshev coefficients, and
leaves the same compact-open limit because $s_N/N^2\to0$.  On each fixed
disk the added factor is $\exp(O(N^{-1/2}))$, so convergence is slower but
the limit is unchanged.  Thus measure and moment convergence do not determine
the microscopic root configuration.
\end{remark}

\section{Scope and limitations}

The results are obstruction theorems. They do not assert existence for a
prescribed target. The independently model-rechecked construction in the Proof
Idea demonstrates that the joint hypotheses are nonempty, but it is not a
general existence theorem. The boundary examples isolate failed hypotheses.
Nor do the results imply:
\begin{itemize}
  \item coefficientwise convergence of the $a_{m,k}$;
  \item reconstruction of the lattice from the weak limit;
  \item convergence of individual roots outside fixed compact windows;
  \item a rate in $k$ for $r_kd_k^2\to0$;
  \item the same conclusions for signed or complex measures.
\end{itemize}
The order of limits in \eqref{eq:budget} is essential: first
$k\to\infty$ for fixed $A$, and only then $A\to\infty$.  No diagonal choice
$A=A(k)$ is justified without additional uniform information.

\section*{AI Use and Provenance Statement}

The original theorem formulations and candidate proofs were generated and
iteratively revised with substantial assistance from
\mbox{Claude Opus 5} (\emph{max} effort, July 2026).  A separate model,
\mbox{GPT-5.6-sol}, was used for
adversarial review, a separate reconstruction of the arguments, boundary
testing, literature triage, and preparation of this manuscript
(\emph{max} reasoning effort, July 2026).  Two external model-audit rounds
were then conducted by five independent adversarial agents coordinated
through \mbox{Claude Fable 5} (July 29, 2026).  They worked from the manuscript
and its released verification materials, without access to the author's
working notes.  The second round reported no mathematical blockers and
identified process, statement, citation, reproducibility, and editorial
corrections.  Those
corrections are incorporated in this release.  An English audit summary is
included in the release archive and recorded in the release metadata and
workflow disclosure.

In August 2026, \mbox{GPT-5.5} at \emph{xhigh} reasoning effort independently
rechecked the focused non-vacuity construction against every hypothesis of the
forced-window theorem. That recheck is a model audit, not human peer review.

Vasily Stodolsky directed the workflow, set its scope and release criteria,
curated the mathematical and verification record, made the final decisions on
claims and revisions, approved this version, and accepts responsibility for its
contents. The mathematical development was materially assisted by AI systems,
whose roles are disclosed in the paper. No independent human peer review is
claimed. The separate-model review is not human peer review, and no AI system
is listed as an author. No future human
review or later version is claimed or scheduled as a condition of releasing
version 0.1.4.

\begingroup
\small
\bibliographystyle{alpha}
\bibliography{references}
\endgroup

\end{document}